\documentclass{article}
\usepackage[T2A]{fontenc}
\usepackage[utf8]{inputenc}
\usepackage[english, russian]{babel}
\usepackage{latexsym}
\usepackage{amssymb}
\usepackage{amsmath}
\usepackage{amsthm}
\usepackage{tikz, subfig}
\usetikzlibrary{decorations.markings}
\usepackage{hyperref}
\usepackage{graphicx}
\usepackage{caption}
\usepackage{leading}
\usepackage[a4paper,bindingoffset=0.2in,%
            left=1.25in,right=1.25in,top=1in,bottom=1.2in,%
            footskip=.25in]{geometry}
\newtheorem*{my_theorem*}{Theorem}
\newtheorem{my_lemma}{Lemma}

\newtheorem{my_def}{Definition}
\newtheorem{my_corollary}{Corollary}
{}

\begin{document}
\tikzset{->-/.style={decoration={
  markings,
  mark=at position #1 with {\arrow{>}}},postaction={decorate}}}
\tikzset{-<-/.style={decoration={
  markings,
  mark=at position (1 - #1) with {\arrow{<}}},postaction={decorate}}}

{\LARGE\bfseries\center Asphericity of~groups defined by graphs\par}

\vspace{0.3cm}

{\scshape\large\center Vadim Bereznyuk \par}

\vspace{0.6cm}
\noindent\makebox[\textwidth][c]{%
\begin{minipage}{10cm}
{\footnotesize A~graph $\Gamma$ labelled by a~set $S$ defines a~group $G(\Gamma)$ whose generators are the~set of~labels $S$ and whose relations are all words which can be read on~closed paths of~this graph. We introduce the~notion of~aspherical graph and prove that such a~graph defines an~aspherical group presentation. This result generalizes a~theorem of~Dominik Gruber on~graphs satisfying graphical $C(6)$-condition and also allows to get new graphical conditions of asphericity analogous to some classical conditions.}
\end{minipage}}
\vspace{0.7cm}

\section{Introduction}
\label{introduction}
\leading{13pt}

Every group can be defined by a~set of~generators and a~set of~relations among these generators. The~latter set can be excessive: there can be some non-trivial identities among its elements. For example, in~a~group $\langle a, b, c \mid ab^{-1}, bc^{-1}, ac^{-1} \rangle$ a~relation $ac^{-1}$ can be written as $ab^{-1}bc^{-1}$. Roughly speaking, a~presentation is called aspherical if all identities among its relations are trivial. It can be formalized in~various ways, so there are quite a~few different definitions of~asphericity (see, for example, \cite{CCH81}).

It is well known that asphericity follows from the~classical small cancellation conditions. In~$2003$ Mikhail Gromov briefly introduced a~graphical analogue of~small cancellation theory in~his paper \cite{Gro03}. After that Yann Ollivier gave a~combinatorial proof of~a~theorem of~Gromov, which in~particular states asphericity of~groups defined by~graphs satisfying graphical $C'(1/6)$-condition \cite{Oll06}. In~$2015$ Dominik Gruber proved asphericity of~groups defined by~graphs satisfy\-ing graphical $C(6)$-condition \cite{Gru15}.

We introduce a~notion of~aspherical graph and suggest to consider it as graphical analogue of diagrammatic asphericity. That notion allows to transfer known classical conditions which imply diagrammatic asphericity to graphical case. We show that not only graphical analogue of condition $C(6)$ implies asphericity of a group but also graphical analogues of conditions ${C(4)\&T(4)}$ and $C(3)\&T(6)$. Moreover we show how a car-crash lemma from  \cite{Kl93} can be applied to prove asphericity in graphical case.

Classical small cancellation theory operates with presentations where every two distinct relations have quite short common parts. In~graphical small cancellation theory, a~group is defined by a~labelled graph. The~set of~generators is the~set of~labels and the~set of~relations is the~set of~all words which can be read on~closed paths of~the~graph. Thus every relation corresponds to a~closed path where this relation can be read. Unlike classical case, two distinct relations can have a~long common part, but only if this common part originates from the~graph. It means that paths of~the~graph corresponding to these relations have the~same common part as relations themselves.

Recall that a~reduction pair in~a~diagram is a~pair of~distinct faces of~the~diagram such that their boundary cycles share a~common edge and such that their boundary cycles, read starting from that edge, clockwise for one of~the~faces and counter-clockwise for the~other, are equal as words. A~spherical diagram is reduced if there are no reduction pairs. If there exists no reduced spherical diagram over a~presentation, then the~presentation is called diagrammatically aspherical. It is well known that presentations satisfying classical small cancellation conditions are diagrammatically aspherical. 

In graphical case, we call a~pair of~faces a~graphical reduction pair if these faces share an~edge originating from the~graph. A~spherical diagram is graphically reduced if there exists no graphical reduction pair in~this diagram. If there exists no graphically reduced diagram over a~presentation whose set of~relations is the~set of~labels of~all simple closed paths of~the~graph, then we call this graph aspherical. It is easy to show that asphericity of~a~graph follows from the~graphical small cancellation $C(6)$-condition.

It turns out that asphericity of~a~graph implies topological asphericity of~the~corresponding group. Thus asphericity of a graph can be considered as a graphical analogue of diagrammatic asphericity that allows to transfer different conditions which imply diagrammatic asphericity to graphical case.

Note that, when we define a~group by a~graph, we can restrict the~set of~relations to the~set of~labels of~all simple closed paths. It does not change the~group. This set of~relations can be reduced further. We can choose an~arbitrary basis of~the~fundamental group of~the~graph. Then the~set of~relations will be the~set of~cyclically reduced paths of~that basis. Again it does not change the~group. The~obtained presentation is topologically aspherical if the~graph is aspherical.

The paper begins with a~brief introduction to theory of~groups defined by graphs. In~Section \ref{main_result} the~main result is formulated and in~Section \ref{idea_of_proof} the~proof is outlined. In~Section \ref{main_notions} we give exact definitions of~main notions. Section \ref{relations_and_spherical_diagrams} is devoted to the~link between identities among relations of~a~presentation and spherical diagrams over this presentation. The~full proof of~the~main theorem can be found in~Section \ref{proof}. At the end we show how to transfer classical conditions of asphericity to graphical case.

The author thanks Anton Klyachko for many useful conversations and remarks, and the anonymous referee for valuable remarks that improved this work.

\subsection{Groups defined by graphs and graphical small \\
cancellation conditions}
\label{graphical_groups_intro}

Let $\Gamma$ be an~oriented graph every edge of~which is labelled by an~element of~a~finite set $S$. Then each path $p$ in that graph can be mapped to a word $\ell(p)$ in the alphabet $S\sqcup S^{-1}$, which is called the~label of~the~path $p$. This word is equal to a product (without reductions) of the labels of edges of this path, considering that if orientation of the edge in the path doesn't match orientation of the edge in the graph then the label belongs to the product with exponent $-1$.

Let $R_c$ be a~set of~labels of~all closed paths in~$\Gamma$, $R_s$ be a~set of~labels of~all simple closed paths in~$\Gamma$ and $R_f$ be a~set of~cyclically reduced labels of~paths which generate a~basis of~a~fundamental group of~each connected component of~the~graph $\Gamma$ (note that $R_c$ and $R_s$ are determined by the~graph $\Gamma$ itself while $R_f$ depends on~the~chosen basis of~the~fundamental group of~$\Gamma$). Then a~group $G(\Gamma)$ are defined by one of~the~three following presentations: $\langle S \mid R_{c} \rangle$, $\langle S \mid R_{s} \rangle$ or $\langle S \mid R_{f} \rangle$. Clearly, all these presentations define the~same group.

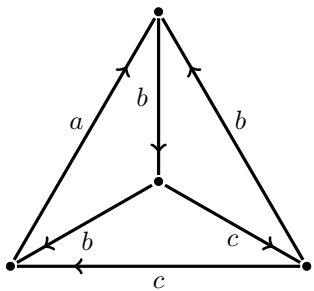
\begin{figure}[!h]
\centering
\begin{tikzpicture}[scale=0.75]

\node[inner sep=0pt, minimum size=4.5pt] (O) at (0, 0) {};
\node[inner sep=0pt, minimum size=4.5pt] (A) at (0, 3) {};
\node[inner sep=0pt, minimum size=4.5pt] (B) at (-0.866 * 3, -1.5) {};
\node[inner sep=0pt, minimum size=4.5pt] (C) at (0.866 * 3, -1.5) {};

\filldraw (O) circle (2pt);
\filldraw (A) circle (2pt);
\filldraw (B) circle (2pt);
\filldraw (C) circle (2pt);

\draw[-<-=0.8, very thick] (O) -- node[left] {$b$} (A);
\draw[->-=0.8, very thick] (O) -- node[below, yshift=1pt, xshift=1pt] {$b$} (B);
\draw[->-=0.8, very thick] (O) -- node[below] {$c$} (C);
\draw[->-=0.8, very thick] (B) -- node[above, xshift=-3pt] {$a$} (A);
\draw[->-=0.8, very thick] (C) -- node[below] {$c$} (B);
\draw[->-=0.8, very thick] (C) -- node[above, xshift=3pt] {$b$} (A);
\end{tikzpicture}

\caption{$G(\Gamma) \cong \langle a, b, c \mid bbc, c^{-1}bc^{-1}, b^{-1}a^{-1}b^{-1} \rangle$}
\label{example_1:fig}
\end{figure}

\begin{my_def}
A lift of~a~word $w$ in~the~graph $\Gamma$ is such a~path $\bar p$ in~the~graph that $\ell(\bar p) \equiv w$ (i.e., the~label of~the~path $\bar p$ coincides with the~word $w$ character by character).
\end{my_def}

\begin{my_def}
A word $w$ is a~piece (with respect to $\Gamma$) if it has two (or more) distinct lifts in~the~graph $\Gamma$.
\end{my_def}

\begin{my_def}
Let $p$ be a~path in~a~graph labelled by a~set $S$. A~lift of~the~path $p$ in~the~graph $\Gamma$ is such a~path $\bar p$ in~the~graph that $\ell(\bar p) \equiv \ell(p)$ (i.e., the~label of~the~path $p$ coincides with the~label of~the~path $\bar p$ character by character).
\end{my_def}

\begin{my_def}
Let $p$ be a~path in~a~graph labelled by a~set $S$. The~path $p$ is a~piece (with respect to $\Gamma$) if it has two (or more) distinct lifts in~the~graph $\Gamma$.
\end{my_def}

Recall that a cycle in a graph is a set of all cyclic shifts of some closed path.

\begin{my_def}
    Let $\gamma$ be a cycle in a graph labeled by a set $S$. A lift of the cycle $\gamma$ in the graph $\Gamma$ is such a cycle $\bar \gamma$ in the graph together with a map $f\colon \gamma \to \bar \gamma$, that $f$ commutes with cyclic shifts and $f(p)$ is a lift of $p$ for all $p \in \gamma$. 
\end{my_def}

Everywhere further it will be clear about which graph $\Gamma$ we talk, so we will call words and paths just “pieces”, not “pieces with respect to $\Gamma$”.

Consider an~example. Let $\Gamma$ be a~graph as in~Figure \ref{example_1:fig}. Then the~words $b$, $b^{-1}$, $c$ and $c^{-1}$ are all pieces of~lengths $1$. The~words $a$ and $a^{-1}$ are not pieces. The~words $bb$ and $(bb)^{-1}$ are all pieces among reduced words of~length $2$. 

A labelling of~a~graph $\Gamma$ is reduced if any two distinct edges starting at the~same vertex have distinct labels and any two distinct edges ending at the~same vertex have distinct labels.

\begin{my_def}
Let $\Gamma$ be a~labelled graph and let $k \in \mathbb{N}$. We say $\Gamma$ satisfies graphical $C(k)$-condition (or $\Gamma$ is a~$C(k)$-graph) if:
\begin{itemize}
\item the~labelling of~$\Gamma$ is reduced and
\item no simple closed path is a~concatenation of~strictly fewer than $k$ pieces.
\end{itemize}
\end{my_def}

Note that if a graph $\Gamma$ satisfies graphical condition $C(2)$ then that graph has a reduced labelling and any word from $R_s$ has a unique lift in the graph.

A graph as in~Figure \ref{example_1:fig} satisfies graphical $C(2)$-condition, but does not satisfy graphical $C(3)$-condition because the~simple closed path with the~label $bbc$ is a~concatenation of~the~pieces $bb$ and $c$.

Let $\Gamma$ be a $C(2)$-graph and let $D$ be a diagram over the presentation $\langle S \mid R_s \rangle$ (see the~next section for definitions). Let $p$ be a path lying in intersection of some positively oriented boundary path of a face $\Pi_1$ and some negatively oriented boundary path of a face $\Pi_2$. A word from $R_s$ are written on the boundary of any face of $D$. Thus the boundary of any face has a lift in the graph $\Gamma$. That lift is unique since $\Gamma$ is a $C(2)$-graph.

Lift the~boundary of~the~face $\Pi_1$ in~the~graph. After that the~path $p$, as subpath of~the boundary, maps to some path $p_1$ in~the~graph. Similarly lift the~boundary of~the~face $\Pi_2$ and determine a~path $p_2$. We say that the~path $p$ originates from the~graph $\Gamma$ if $p_1 = p_2$. Roughly speaking, a~path $p$ originates from the~graph $\Gamma$ if faces $\Pi_1$ and $\Pi_2$ share the~same path in~the~diagram and in~the~graph $\Gamma$ itself.

Note that if a~path $p$, lying between faces $\Pi_1$ and $\Pi_2$, does not originate from the~graph $\Gamma$ then it is a~piece. Indeed, if $p$ does not originate from the~graph then its lifts via $\Pi_1$ and via $\Pi_2$ are distinct. Therefore this path have two distinct lifts in~the~graph $\Gamma$, i.e., this path is a~piece.

\begin{figure}[!t]
\centering

\subfloat{
\begin{tikzpicture}[scale=0.66, baseline]

\node[inner sep=0pt, minimum size=4.5pt] (O) at (0, 0) {};
\node[inner sep=0pt, minimum size=4.5pt] (A) at (0, 3) {};
\node[inner sep=0pt, minimum size=4.5pt] (B) at (-0.866 * 3, -1.5) {};
\node[inner sep=0pt, minimum size=4.5pt] (C) at (0.866 * 3, -1.5) {};
\node () at (-3.5, 2.5) {\Large $\Gamma$};

\filldraw (O) circle (2pt);
\filldraw (A) circle (2pt);
\filldraw (B) circle (2pt);
\filldraw (C) circle (2pt);

\draw[-<-=0.8, very thick] (O) -- node[left] {$b$} (A);
\draw[->-=0.8, very thick] (O) -- node[below, yshift=2pt, xshift=2pt] {$b$} (B);
\draw[->-=0.8, very thick] (O) -- node[below] {$c$} (C);
\draw[->-=0.8, very thick] (B) -- node[above, xshift=-3pt] {$a$} (A);
\draw[->-=0.8, very thick] (C) -- node[below] {$c$} (B);
\draw[->-=0.8, very thick] (C) -- node[above, xshift=3pt] {$b$} (A);
\end{tikzpicture}
}
\quad\quad\quad\quad
\subfloat{
\begin{tikzpicture}[scale=0.66, baseline]

\node[inner sep=0pt, minimum size=4.5pt] (O) at (0, 0) {};
\node[inner sep=0pt, minimum size=4.5pt] (A) at (0, 3) {};
\node[inner sep=0pt, minimum size=4.5pt] (B) at (-0.866 * 3, -1.5) {};
\node[inner sep=0pt, minimum size=4.5pt] (C) at (0.866 * 3, -1.5) {};
\node[inner sep=0pt, minimum size=4.5pt] (D) at (0.866 * 3, 1.5) {};
\node[inner sep=0pt, minimum size=4.5pt] (E) at (-0.866 * 3, -3) {};
\node[inner sep=0pt, minimum size=4.5pt] (F) at (0.866 * 3, -3) {};
\node () at (-3.5, 2.5) {\Large $D$};

\filldraw (O) circle (2pt);
\filldraw (A) circle (2pt);
\filldraw (B) circle (2pt);
\filldraw (C) circle (2pt);
\filldraw (D) circle (2pt);
\filldraw (E) circle (2pt);
\filldraw (F) circle (2pt);

\draw[-<-=0.8, very thick, dashed] (O) -- node[left] {$b$} (A);
\draw[->-=0.8, very thick, dashed] (O) -- node[below, yshift=2pt, xshift=2pt] {$b$} (B);
\draw[->-=0.8, very thick, dashed] (O) -- node[below] {$c$} (C);
\draw[->-=0.8, very thick] (B) -- node[above, xshift=-3pt] {$a$} (A);
\draw[->-=0.8, very thick] (C) -- node[below] {$c$} (B);
\draw[->-=0.8, very thick] (C) -- node[above, xshift=3pt] {$b$} (A);
\draw[->-=0.8, very thick] (A) -- node[above] {$b$} (D);
\draw[->-=0.8, very thick] (D) -- node[right] {$a$} (C);
\draw[->-=0.8, very thick] (B) -- node[left] {$b$} (E);
\draw[-<-=0.8, very thick] (E) -- node[below] {$a$} (F);
\draw[-<-=0.8, very thick] (F) -- node[right] {$b$} (C);
\end{tikzpicture}
}

\caption{A graph $\Gamma$ and a~diagram $D$ over $\langle a, b, c \mid R_s \rangle$. Dotted edges originate from the~graph.}
\end{figure}
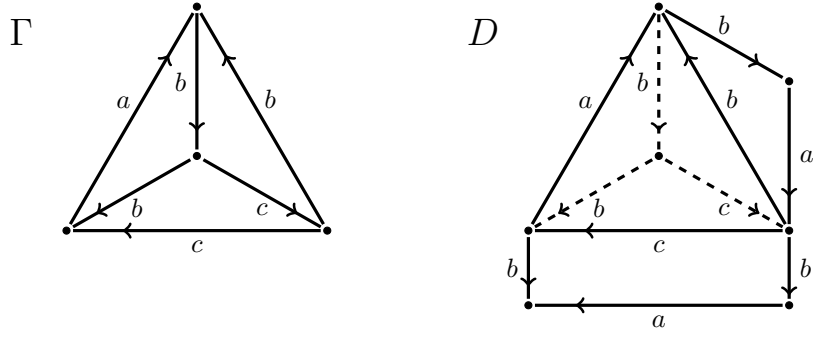

\subsection{Main result}
\label{main_result}

Recall that a~presentation complex $K(S; R)$ of~a~presentation $\langle S \mid R \rangle$ is a~$2$-complex which has a~$1$-skeleton which consists of~a~single vertex and a~loop labelled by $s$ for every element $s$ from $S$, and which have a~face with the~boundary label $r$ attached to the~$1$-skeleton for every element $r$ from $R$. 

\begin{my_def}
A presentation $\langle S \mid R \rangle$ is called aspherical if its presentation complex $K(S; R)$ is aspherical, i.e., $\pi_q(K(S; R)) = 0$ for all $q \geq 2$.
\end{my_def}

\begin{my_def}
We say that a~diagram $D$ over the~presentation $\langle S \mid R_s \rangle$ is graphically reduced if it does not have edges originating from the~graph $\Gamma$.
\end{my_def}

\begin{my_def}
We say that a~graph $\Gamma$ is aspherical if it satisfies graphical $C(2)$-condition and there exists no graphically reduced spherical diagram over the~presentation $\langle S\mid R_s \rangle$.
\end{my_def}

The following theorem is the~main result of~this paper.

\begin{my_theorem*}
If a~graph $\Gamma$ is aspherical, then the~presentation $\langle S\mid R_f \rangle$ is aspherical.
\end{my_theorem*}

\subsection{Idea of~proof}
\label{idea_of_proof}

Propositions $1.3$ and $1.5$ of~the~paper \cite{CCH81} implies that a~presentation $\langle S \mid R \rangle$, where all relations from $R$ are not empty and freely reduced, is aspherical if and only if the~presentation is concise, no relation is a~proper power and any identity among relations of~this presentation is trivial.

Recall that a~presentation $\langle S \mid R \rangle$ is concise if for any two distinct relations $r$ and $r'$ from $R$ nor $r$, neither $r^{-1}$ is conjugate to $r'$.

Also recall a~notion of~identity among relations of~a~presentation $\langle S \mid R \rangle$. Let $\pi = (p_1, \dots, p_n)$ be a~sequence such that $p_i = u_i r_i^{\epsilon_i} u_i^{-1}$, where $r_i \in R$, $u_i \in F(S)$ and $\epsilon_i \in \{+1, -1\}$. It is called an~identity if a~product of~its elements is equal to the~identity element of~the~free group, i.e., $p_1 \cdots p_n = 1$ in~$F(s)$. There are identities which we should consider as trivial. For this reason Peiffer transformations are introduced:
\begin{enumerate}
\item Replace any pair of~consecutive elements $(p_i, p_{i+1})$, either by the~pair $(p_i p_{i+1} p_i^{-1}, p_i)$ or by the~pair $(p_{i+1}, p_{i+1}^{-1} p_i p_{i+1})$.
\item Delete the~pair of~consecutive elements $(p_i, p_{i+1})$ if $p_i p_{i+1} = 1$ in~$F(S)$.
\item Insert at any place a~pair of~inverse elements $(p, p^{-1})$.
\end{enumerate}

An identity is called trivial if it can be transformed to an~empty identity by the~finite number of~Peiffer transformation.

Lemma $2.22$ of~the~paper \cite{Gru15} guarantees the~first two conditions: conciseness of~a~presentation and absence of~proper powers. Thus we only should show that any identity among relations is trivial.

To prove this fact we use a~link between identities among relations of~a~presentation and spherical diagrams over this presentation which was obtained in~\cite{CH82}. A~plan of~the~proof is the~following. Assume the~contrary, that there exists non-trivial identities over the~presentation $\langle S\mid R_f \rangle$ or, equivalently, that there exists non-trivial spherical diagrams over the~presentation $\langle S\mid R_f \rangle$. Consider a~part of~an~$1$-skeleton of~a~spherical diagram which consists of~all edges not originating from a~graph $\Gamma$. We call such a~part a~not originating skeleton. Consider a~non-trivial spherical diagram with the~smallest not originating skeleton. Delete all edges originating from the~graph from this diagram. It turns out that the~obtained diagram is a~diagram over the~presentation $\langle S\mid R_s \rangle$. Moreover, it is graphically reduced because we deleted all originating edges. A~contradiction with asphericity of~the~graph $\Gamma$.

\section{Main notions}
\label{main_notions}

Definitions of~graphs and diagrams are given in~this section.

\subsection{Graphs}
\label{graphs}

We use a~definition of~graph according to \cite{LiShu77}, i.e., a~graph is a~union of~two sets $V$ and $E$ with three maps $\alpha \colon E \to V$, $\omega \colon E \to V$, $\cdot^{-1}\colon E \to E$. The~elements of~$V$ are called vertices and the~elements of~$E$ are called edges. If $e \in E$ then $\alpha(e)$ is called the~initial point and $\omega(e)$ the~terminal point. The~map $\cdot^{-1}$ assigns to every edge its inverse. For convenience we write $\cdot^{-1}(e)$ as $e^{-1}$. The~map $\cdot^{-1}$ should be an~involution without fixed elements such as $\alpha(e^{-1}) = \omega(e)$ and $\omega(e^{-1}) = \alpha(e)$. In~fact, it is a~definition of~undirected graph, because for every edge $e$ graph also contains $e^{-1}$.

A path in a graph is a finite sequence of edges $p = (e_1, \dots, e_n)$, such that $\alpha(e_{i+1}) = \omega(e_i)$ for $1 \leqslant i < n$. A path begins at a point $\alpha(p) = \alpha(e_1)$ and ends at a point $\omega(p) = \omega(e_n)$. 

A labelling of~a~graph $\Gamma$ by a~set $S$ is such a~map $\ell\colon E \to S\sqcup S^{-1}$ that $\ell(e^{-1}) = {\ell(e)}^{-1}$. A~labelled graph is a~graph with its labelling. One can think of~a~labelled graph as an~oriented graph every edge of~which labelled by an~element of~$S$. A~labelling of~a~graph $\Gamma$ is reduced if any two distinct edges starting at the~same vertex have distinct labels and any two distinct edges ending at the~same vertex have distinct labels.

Continue the~map $\ell$ to a~set of~all paths in~the~graph $\Gamma$. Let $p = (e_1, \dots, e_n)$ be a~path in~the~graph, then put $\ell(p) = \ell(e_1) \cdots \ell(e_n)$, where $w_1 \cdot w_2$ is concatenation. Thus $\ell(p)$ is a~word over an~alphabet $S \sqcup S^{-1}$ (not necessarily reduced). We call $\ell(p)$ as the~label of~the~path $p$. 

Let $p = (e_1, \dots, e_n)$ be a~path in~a~graph. The~path $p$ is reduced if it contains no subpaths $(e, e^{-1})$. The~path is trivial if it becomes empty after consecutive deletion of~all subpaths $(e, e^{-1})$. The~path is closed if its initial point coincides with its terminal point or if it is empty. The~path is simple if it is not empty and it does not contain non-empty closed subpaths. The~path $p$ is simple closed if it is not trivial, is closed and no proper subpaths of~$p$ is closed. A set of cyclic shifts of a closed path is called a cycle. A path $p$ is called an arc if all its vertices besides endpoints have degree $2$. An arc is called a spur if at least one its endpoint has degree $1$.

\subsection{Diagrams over graphs}
\label{diagrams}

Let $R$ be a set of words, then we define $R_{sym}$ as the~set obtained from the~sets $R$ and $R^{-1}$ by considering all its elements up to cyclic shifts.

A singular disk diagram in~the~alphabet $S$ is a~finite and simply connected $2$-dimensional CW-complex embedded into $\mathbb{R}^2$ such that its $1$-skeleton is a~graph labelled by the~set $S$. The~closures of~its $1$-cells and $2$-cells are called edges and faces, respectively. The~label of~a~face is a~cyclic word obtained by reading its boundary path in~a~counterclockwise direction. A~singular disk diagram over a~presentation $\langle S \mid R \rangle$ is a~singular disk diagram in~the~alphabet $S$ every face of~which has a~label from $R_{sym}$. A~simple disk diagram is a~singular disk diagram homeomorphic to a~disk.

We introduce spherical diagrams following \cite{CH82}. First we inductively define a~spherical complex. A~$2$-sphere is a~spherical complex. If $D$ is a~spherical complex, then $D'$ obtained by~attaching a~simple curve or a~$2$-sphere to a~point of $D$ is also a~spherical complex. A~spherical complex is a~complex which can be obtained from a~$2$-sphere by the finite number of~such attachments. In~other words, a~spherical complex is a~tree embedded into $\mathbb{R}^3$ some vertices of~which (at least one) are replaced by spheres or by some number of~spheres attached to each other. A~spherical diagram over a~presentation $\langle S \mid R \rangle$ is a~$2$-dimensional CW-complex homeomorphic to some spherical complex, such that its $1$-skeleton is a~graph labelled by the~set $S$ and such that the~labels of~all its faces lie in~$R_{sym}$. A~simple spherical diagram is a~spherical diagram homeomorphic to a~sphere 

Let denote a cycle consisting of~positively oriented (i.e., obtained by reading the~boundary in~a~counterclockwise direction) boundary paths of~a~face $\Pi$ as $\partial\Pi^{+}$, and a cycle consisting of~negatively oriented as $\partial\Pi^{-}$. It is clear that $\partial\Pi^{-} = \{\gamma^{-1}\colon\gamma \in \partial\Pi^{+}\}$. If $P$ and $Q$ are two sets of~paths then $P \Cap Q = \{r \colon \exists p \in P, \exists q \in Q \text{ such as } p = r p', q = r q'\}$.

\begin{my_def}
Let $D$ be a~diagram over the~presentation $\langle S\mid R_c\rangle$. Let $\Pi_1$ and $\Pi_2$ be two faces of~$D$ (not necessary distinct)  and let $p \in \partial\Pi_1^{+} \Cap \partial\Pi_2^{-}$. We say that the~path $p$ originates from the~graph $\Gamma$ if there exists such lifts of~$\partial\Pi_1^{+}$ and $\partial\Pi_2^{-}$ in~$\Gamma$, that a~lift of~$p$ in~$\Gamma$ via $\partial\Pi_1^{+}$ and via $\partial\Pi_2^{-}$ are equal. 
\end{my_def}

Note that according to this definition edges on the boundary of a diagram are not originating. And note again that if the~graph $\Gamma$ satisfies graphical $C(2)$-condition then lifts of~$\partial\Pi_1^{+}$ and $\partial\Pi_2^{-}$ are unique.

\section{Identities and spherical diagrams}
\label{relations_and_spherical_diagrams}

For any identity $\pi = (p_1,\dots,p_n)$ over a~presentation $\langle S \mid R \rangle$ we can construct a~spherical diagram over the~same presentation by the~so-called van Kampen construction (\cite{LiShu77}, \cite{CH82}). Recall it. First, note that for any word $w$ in the alphabet $S\sqcup S^{-1}$ a linear graph $p$ labelled by the set $S$ such that $\ell(p) \equiv w$ can be constructed.

Now, fix a~point $v$ on~the~plane. After that for each $p_i = u_i r_i^{\varepsilon_i} u_i^{-1}$ draw a face with a spur on the plane such that boundary of the face is a closed arc with a label $r_i^{\varepsilon_i}$ and the spur has a label $u_i$. Then the~boundary label of~the~obtained diagram is equal to the~word $p_1 \cdots p_n$. So this boundary label is trivial because $\pi$ is an~identity. It means that there exists $2$ consecutive edges with opposite labels on~the~boundary. Glue these edges together. We obtain again a~diagram with a~trivial boundary label but with the~smaller boundary size. Consecutively gluing pairs of~edges with opposite labels we finally obtain a~spherical diagram $D$ which is called a~diagram for $\pi$ (see \cite[Section 1.5]{CH82} for details).

Note that obtained diagram is spherical by our definition because if a sphere arises after gluing then this sphere touches the boundary of a diagram only by one vertex.

This procedure can be reversed. By “ungluing” faces of~a~spherical diagram $D$ over a~presentation $\langle S \mid R \rangle$ along edges we can obtain a~bouquet of~faces with spurs that gives us a~sequence of~elements $\pi$ over the~same presentation. Moreover, $\pi$ is an~identity because the~diagram was spherical. Note that $D$ is a~diagram for $\pi$ because we can reverse the~described procedure.

These procedures are not unique. We may obtain different diagrams and identities changing the~order of~gluings and ungluings. But at the~same time Proposition $8$ from \cite{CH82} implies that if $D$ is a~diagram for $\pi_1$ and a~diagram for $\pi_2$ then $\pi_1$ is trivial if and only if $\pi_2$ is. It allows us to define trivial spherical diagrams with correct correspondence to trivial identities. We will call $D$ a~trivial spherical diagram if a~trivial identity can be obtained from $D$. Due to the~fact noted above, every identity obtained from a~trivial spherical diagram $D$ will be trivial. Thus only trivial identities can be obtained from a~trivial diagram, and only trivial diagrams can be obtained from a~trivial identity.

Note also that actually a~diagram $D$ which is a~diagram for some identity $\pi$ has the~marked point (the initial vertex $v$ from the~construction). Intuitively it is clear that replacing of the marked point does not change triviality of a diagram. To strictly prove that we should use original definition of triviality of diagram from \cite{CH82}, which states that a trivial diagram is a diagram which can be transformed to the trivial diagram consisting of only one vertex by finite number of certain transformations. Note that these transformations do not depend on the marked point that means that if a diagram can be transformed to trivial then after replacing of the marked point it still can be transformed to trivial. Thus replacing of the marked point does not change triviality of a diagram so further we will not specify which point we consider as marked.

And finally note that a~spherical diagram can have spurs. But Proposition $8$ from \cite{CH82} implies that inserting and deleting of~spurs does not change triviality of~the~diagram.

\section{Proof of~the~theorem}
\label{proof}

In this section $\Gamma$ is a~labelled by the~set $S$ graph and $R_c$, $R_s$ and $R_f$ are the~sets defined at the~beginning of~the~section \ref{graphical_groups_intro}.

\begin{my_def}
Let $D$ be a~diagram over the~presentation $\langle S\mid R_s \rangle$. A~not originating skeleton of~$D$ is a~graph that consists of~all edges of~$D$ which do not originate from the~graph $\Gamma$.
\end{my_def}

\begin{my_lemma}
Suppose $\Gamma$ satisfies graphical $C(2)$-condition and let $D$ be a~simple spherical diagram over the~presentation $\langle S\mid R_s \rangle$. Then a~not originating skeleton of~$D$ has no spurs.
\end{my_lemma}

\begin{proof}[Proof]
Assume the~contrary. Let $D$ be a~simple spherical diagram whose not originating skeleton contains a~spur $e$ with endpoint $v$. It means that all other edges incidental to $v$ originate from the~graph $\Gamma$ because otherwise the~edge $e$ would not be a~spur in~the~originating skeleton.

We may assume that a~part of~the~diagram in~a~neighborhood of~the~point $v$ is embedded into the~plane. Draw a~circle with a~center at $v$ with quite a~small radius such that the~circle intersects all incidental to $v$ edges and only them. Once if an~edge is not a~loop and twice if an~edge is a~loop. This circle intersects some number of consecutive faces $\Pi_1, \dots, \Pi_n$. Also it intersects the edges $e_1, \dots, e_{n-1}, e_n = e$, where $e_i$ lies in~the~intersection of~faces $\Pi_i$ and $\Pi_{i+1}$.

The boundary of~each face has a~unique lift to the~graph because $\Gamma$ satisfies graphical $C(2)$-condition. Let $v_i$ be a~vertex of~the~graph obtained by lifting $v$ via the~boundary of~the~face $\Pi_i$. The~edge $e_i$ originates from the~graph for $i = 1, \dots, n-1$, therefore $v_i = v_{i+1}$ for $i = 1, \dots, n-1$. Thus $v_1 = v_n$. It means that lifts of~the~edge $e$ via $\partial \Pi_1^{+}$ and via $\partial \Pi_1^{+}$ have a~common vertex in~the~graph and, moreover, have the~same label. It implies that these lifts coincide because a~labelling of~the~graph is reduced. Thus $e$ originates from the~graph. A~contradiction.
\end{proof}

\begin{figure}[!t]
\centering
\begin{tikzpicture}
[node_style/.style={shape=circle,draw,fill=black!100,inner sep=1pt, minimum size=3pt}]

\node () at (1.1, -1.1) [shape=circle, fill=white!100] {$\Pi_1$};
\node () at (1.1, 1.1) [shape=circle, fill=white!100] {$\Pi_2$};
\node () at (-1.1, 1.1) [shape=circle, fill=white!100] {$\Pi_3$};
\node () at (-1.1, -1.1) [shape=circle, fill=white!100] {$\Pi_4$};
\node () at (0.22, 0.22) [shape=circle, fill=white!100] {$v$};
\node (B) at ( 2,0) [node_style] {};
\node (C) at ( 0,2) [node_style] {};
\node (D) at (-2,0) [node_style] {};
\node (E) at ( 0,-2) [node_style] {};
\node (A) at ( 0,0) [node_style] {}
	edge [-, dashed, very thick] node[auto] {$e_1$} (B)
	edge [-, dashed, very thick] node[auto] {$e_2$} (C)
	edge [-, dashed, very thick] node[auto] {$e_3$} (D)
	edge [-, solid, very thick] node[auto] {$e$} (E);
\draw[->-=0.15, thick] (0, 0) circle [radius=0.65];
\end{tikzpicture}
\caption{Illustration for Lemma $1$. Dotted edges originate from the~graph.}
\end{figure}
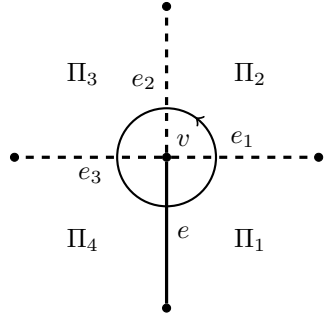

\begin{my_def}
We will call $D$ a~minimal non-trivial spherical diagram over the~presentation $\langle S\mid R_f \rangle$ if it is a~non-trivial spherical diagram such that
\begin{enumerate}
\item A~not originating skeleton of~$D$ has the~smallest number of~edges among all non-trivial diagrams and among all such diagrams a~not originating skeleton of~$D$ has the~smallest number of~vertices.
\item $D$ has the~smallest number of~edges among all diagrams satisfying the~first condition.
\end{enumerate}
\end{my_def}

Note that a~minimal non-trivial spherical diagram $D$ is always a~simple spherical diagram. Indeed, assume the~contrary, let $D$ not be a~simple spherical diagram. We noted earlier that insertion and deletion of~spurs does not affect its triviality. Therefore $D$ does not have spurs. It means that there exists such a~point $v$ in~$D$ that the~diagram splits into two spherical diagrams $D_1$ and $D_2$ after cutting $D$ at $v$. Let $\pi_1$ be an~identity obtained from the~diagram $D_1$ and $\pi_2$ be an~identity obtained from the~diagram $D_2$. These identities are trivial because $D$ is a~minimal non-trivial spherical diagram. Then an~identity $\pi = (\pi_1, \pi_2)$ is trivial as well. But $\pi$ can be obtained from the~diagram $D$. Thus $D$ is a~trivial spherical diagram that contradicts its definition.

Let $D$ be a~simple spherical diagram. Denote by $\bar D$ a~simple spherical diagram obtained from $D$ by erasing all edges originating from the~graph $\Gamma$ (if isolated vertices are left we delete them). Clearly, a~$1$-skeleton of~$\bar D$ coincides with a~not originating skeleton of~$D$. Note that actually $\bar D$ may have some not simply connected “faces”, but we still consider them as faces and call them not simply connected faces. Note that the~boundary of~a~not simply connected face consists of~some connected components.  

\begin{my_lemma}
Let $\Gamma$ be a~$C(2)$-graph. Let $D$ be a~minimal non-trivial spherical diagram over the~presentation $\langle S\mid R_f \rangle$. Then the~boundary label of~any face of~the~diagram $\bar D$ is reduced (the boundary label of~any connected component for not simply connected faces).
\end{my_lemma}

\begin{proof}[Proof]
Assume the~contrary. Let $\bar \Pi$ be a~face with a~not reduced boundary label. Let $e$ and $f$ be the~edges where reduction is occured. We will use diamond moves introduced in~\cite{CH82}.

For definiteness assume that $\omega(e) = \omega(f)$ and $\ell(e) = \ell(f) = a$. Consider in~detail only the~case when all $3$ points $\omega(e), \alpha(e), \alpha(f)$ are distinct. Transform the~diagram $D$ as shown in~Fig. \ref{lemma_3_case_1:fig}. Note that the~not originating skeleton was reduced at least by one edge because one of~the~two new edges appears to be inside the~face $\Bar \Pi$. But the~not originating skeleton does not have spurs so this edge originates from the~graph. 

\begin{figure}[!h]
\centering
\subfloat{
\begin{tikzpicture}
[node_style/.style={shape=circle,draw,fill=black!100,inner sep=1pt, minimum size=3pt}]

\node () at (-1.1, 1.4) [shape=circle, fill=white!100] {$\Bar\Pi$};
\node () at (0, 2) {};
\node () at (0, -1.94) {} ;
\node (B) at ( 1,0) [node_style] {};
\node (C) at ( 0,1) [node_style] {};
\node (D) at (-1,0) [node_style] {};
\node (E) at ( 0,-1) [node_style] {};
\node (F) at ( -2, 0) [node_style] {}
	edge [-, solid, very thick] (D);
\node (G) at ( 2, 0) [node_style] {}
	edge [-, solid, very thick] (B);
\node (A) at ( 0,0) [node_style] {}
	edge [->-=0.8, solid, very thick] node[auto] {$a$} (B)
	edge [-, dashed, very thick] (C)
	edge [->-=0.8, solid, very thick] node[auto, swap] {$a$} (D)
	edge [-, solid, very thick] (E);
\end{tikzpicture}
}
\quad
\subfloat{
\begin{tikzpicture}
[node_style/.style={shape=circle,draw,fill=black!100,inner sep=1pt, minimum size=3pt}]

\node () at (-1.1, 1.4) [shape=circle, fill=white!100] {$\Bar\Pi$};
\node (A) at (0, 1) [node_style] {};
\node (B) at (1, 0) [node_style] {};
\node (C) at (0, -1)  [node_style] {};
\node (D) at (-1, 0) [node_style] {};
\path
	(A) edge[->-=0.8, very thick] node[auto] {$a$} (B)
	(A) edge[->-=0.8, very thick] node[auto, swap] {$a$} (D)
	(C) edge[->-=0.8, very thick] node[auto, swap] {$a$} (B)
	(C) edge[->-=0.8, very thick] node[auto] {$a$} (D);
\node (E) at (0, 2) [node_style] {}
	edge [-, dashed, very thick] (A);
\node (F) at (2, 0) [node_style] {}
	edge [-, solid, very thick] (B);
\node (G) at (0, -2) [node_style] {}
	edge [-, solid, very thick] (C);
\node (G) at (-2, 0) [node_style] {}
	edge [-, solid, very thick] (D);
\end{tikzpicture}
}
\quad
\subfloat{
\begin{tikzpicture}
[node_style/.style={shape=circle,draw,fill=black!100,inner sep=1pt, minimum size=3pt}]

\node () at (-0.6, 1.4) [shape=circle, fill=white!100] {$\Bar\Pi$};
\node (O) at (0, 0) [node_style] {};
\node (A) at (0, 1) [node_style] {};
\node (B) at (1, 0) [node_style] {};
\node (C) at (0, -1)  [node_style] {};
\node (D) at (-1, 0) [node_style] {};
\path
	(A) edge[->-=0.8, dashed, very thick] node[auto, swap] {$a$} (O)
	(O) edge[-, very thick] (B)
	(C) edge[->-=0.8, solid, very thick] node[auto] {$a$} (O)
	(O) edge[-, very thick] (D);
\node (E) at (0, 2) [node_style] {}
	edge [-, dashed, very thick] (A);
\node (F) at (0, -2) [node_style] {}
	edge [-, solid, very thick] (C);
\end{tikzpicture}
}
\caption{The first case of~the~Lemma $2$.}
\label{lemma_3_case_1:fig}
\end{figure}
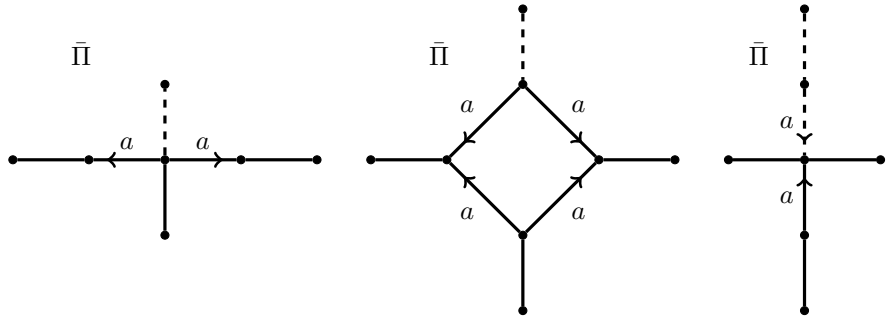

Diamond moves for other cases can be found in~\cite[Section 1.4]{CH82}. We do not consider these cases in~detail because they do not differ much from the~first case: after applying a~diamond move two old edges are replaced by two new ones and one of~the~two new edges appears to be inside the~face $\Bar \Pi$ and therefore originates from the~graph. Thus the~not originating skeleton is reduced anyway.
\end{proof}

\begin{my_lemma}
Let $p$ be a~simple closed path in~the~graph $\Gamma$. Then there exists a~diagram $D$ over the~presentation $\langle S \mid R_f \rangle$ such that $\partial D$ lifts to the~path $p$ and all internal edges of~$D$ originates from the~graph $\Gamma$.
\end{my_lemma}

\begin{proof}[Proof]
Forget for a~while that $\Gamma$ is labelled by the~set $S$. Label every edge $e$ of~the~graph $\Gamma$ by the~new unique label $t_e$. Assume the~path $p$ lies in~the~connected component $\Gamma'$ of~the~graph $\Gamma$ and let $B'$ be the~basis of~$\Pi(\Gamma', v)$ from which some relations of~$R_f$ are obtained. Let $q$ be a~path which runs from the~vertex $\alpha(p)$ to the~vertex $v$ and let $B = \{qb'q^{-1} \mid b' \in B'\}$. Then the~path $p$ are generated by the~paths from $B$, i.e, $p = b_{1}^{\epsilon_1} \cdots b_{n}^{\epsilon_n}$, $b_i \in B$. Each $b_i = u_i r_i u_i^{-1}$ where $r_i$ is a~cyclic reduction of~the~path $b_i$. Due to the~definition of~the~set $R_f$, the~label of~each path $r_i$ lie in~$R_f$.

Fix a~point $s$ on~the~plane. For each $b_{i}^{\epsilon_i}$ draw a~face with a~spur such that the~spur starts at $s$ and has a~label $\ell(u_i)$ and such that the~boundary label of~the~face is $\ell(r_i^{\epsilon_i})$. Thus we obtain a~diagram $E$ such that its boundary label is freely equal to the~label of~the~path $p$. If there exists a~pair of~consecutive edges with the~same labels and opposite directions on~the~boundary of~$E$ then glue them. Doing this several times we obtain a~reduced boundary label. The~labelling of~$\Gamma$ is reduced and $p$ is a~simple closed path so its label is reduced. Thus the~boundary label of~$E$ and the~label of~$p$ is equal as words. 

Note that all internal edges of~$E$ originate from the~graph because they have unique labels and therefore have unique lifts to the~graph. Now recall the~original labelling by the~set $S$. After replacing unique labels by original ones all internal edges still originate from the~graph (because boundaries of~the~faces still have unique lifts to the~graph and these lifts coincide with lifts to the~graph with unique labels) and the~boundary label of~the~diagram still is equal to the~label of~$p$ as words. Thus the~desired diagram for the~path $p$ over the~presentation $\langle S \mid R_f \rangle$ was obtained.
\end{proof}

\begin{my_lemma}
Let $\Gamma$ be a~$C(2)$-graph. Let $D$ be a~minimal non-trivial spherical diagram over the~presentation $\langle S\mid R_f \rangle$ and let $\Pi$ be a~face of~$\bar D$. Then if $\Pi$ is simply connected and its boundary is a~simple closed path then this boundary lifts to a~simple closed path in~the~graph~$\Gamma$.
\end{my_lemma}

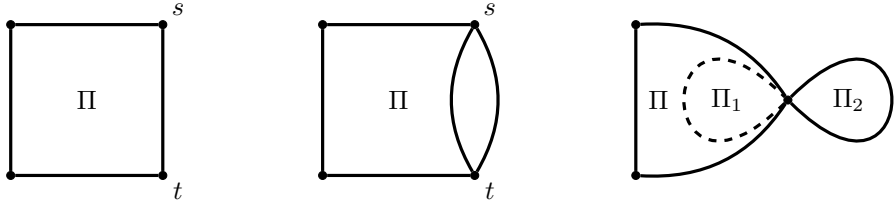
\begin{figure}[!t]
\centering
\subfloat{
\begin{tikzpicture}
[node_style/.style={shape=circle,draw,fill=black!100,inner sep=1pt, minimum size=3pt}]

\node () at (0, 0) {$\Pi$};
\node () at (0, -1.71) {};
\node (A) at (-1, 1) [node_style] {};
\node (B) at (1, 1) [node_style] {};
\node () at (1.2, 1.2) {$s$};
\node (C) at (1, -1) [node_style] {};
\node () at (1.2, -1.2) {$t$};
\node (D) at (-1, -1) [node_style] {};
\path
	(A) edge[-, solid, very thick] (B)
	(B) edge[-, solid, very thick] (C)
	(C) edge[-, solid, very thick] (D)
	(D) edge[-, solid, very thick] (A);
\end{tikzpicture}
}
\qquad\qquad
\subfloat{
\begin{tikzpicture}
[node_style/.style={shape=circle,draw,fill=black!100,inner sep=1pt, minimum size=3pt}]

\node () at (0, 0) {$\Pi$};
\node () at (0, -1.71) {};
\node (A) at (-1, 1) [node_style] {};
\node (B) at (1, 1) [node_style] {};
\node () at (1.2, 1.2) {$s$};
\node (C) at (1, -1) [node_style] {};
\node () at (1.2, -1.2) {$t$};
\node (D) at (-1, -1) [node_style] {};
\path
	(A) edge[-, solid, very thick] (B)
	(B) edge[-, solid, very thick, bend left] (C)
	(B) edge[-, solid, very thick, bend right] (C)
	(C) edge[-, solid, very thick] (D)
	(D) edge[-, solid, very thick] (A);
\end{tikzpicture}
}
\qquad\qquad
\subfloat{
\begin{tikzpicture}
[node_style/.style={shape=circle,draw,fill=black!100,inner sep=1pt, minimum size=3pt}]

\node () at (-0.7, 0) {$\Pi$};
\node (A) at (-1, 1) [node_style] {};
\node (B) at (1, 0) [node_style] {};
\node (D) at (-1, -1) [node_style] {};
\path
	(A) edge[-, solid, very thick, bend left] (B)
	(D) edge[-, solid, very thick, bend right] (B)
	(A) edge[-, solid, very thick] (D)
	(B) edge[-, in=135, out=225, looseness=0.9, distance=2.5cm, dashed, very thick] (B)
	(B) edge[-, in=45, out=-45, looseness=0.9, distance=2.5cm, very thick] (B);
\node () at (0.2, 0) {$\Pi_1$};
\node () at (1.8, 0) {$\Pi_2$};
\end{tikzpicture}
}
\caption{Illustration for Lemma 4.}
\label{lemma_4:fig}
\end{figure}

\begin{proof}[Proof]

Assume the~contrary, that there exists a~face $\Pi$ in~$\bar D$ such that it have a~simple boundary which lifts to a~not simple closed path in~the~graph. By Lemma $2$, the~boundary label of~$\Pi$ is reduced so there exists a~subpath $p$ of~the~path $\partial \Pi$ which lifts to a~simple closed path $\tilde{p}$ in~the~graph. Let $s$ and $t$ be respectively the~start and the~end of~the~path $p$. Transform the~diagram as shown in~Fig. \ref{lemma_4:fig}: cut the~diagram along the~path $p$ and glue together $s$ and $t$. Both $\partial\Pi_1$ and $\partial\Pi_2$ lift to the~simple closed path $\tilde{p}$ in~the~graph $\Gamma$. By the~previous lemma, there exists a~diagram $E$ over the~presentation $\langle S \mid R_f \rangle$ such that its boundary coincides with $\tilde{p}$ and all internal edges originate from the~graph. Glue $E$ at the~place of~$\Pi_1$ and a~diagram symmetric to $E$ at the~place of~$\Pi_2$. Thus we again obtain a~spherical diagram over the~presentation $\langle S \mid R_f \rangle$. Denote it by $D'$.

Note that $\Pi\cap\Pi_1 = \partial\Pi_1$ originates from the~graph because it lifts to the~simple closed path in~the~graph and due to the~graphical $C(2)$-condition such a~lift is unique. Therefore $D'$ has a~smaller not originating skeleton than $D$ because we glued together the~vertices $s$ and $t$ which lay in~a~not originating skeleton of~$D$ and so we reduced the~number of~vertices in~the~not originating skeleton while the~number of~not originating edges are still the~same, because all internal edges of~$E$ originate from the~graph.

Moreover, $D'$ is not trivial because we glued into $D$ the~set of~opposite faces and we can think of~it as adding to the~not trivial identity $\pi$ some sequence $(p_1, \dots, p_n, p_n^{-1}, \dots, p_1^{-1})$ that corresponds to the~$n$-fold application of~the~insert Peiffer transformation. Thus we obtained a~not trivial spherical diagram with a~smaller not originating skeleton than $D$ has. A~contradiction.
\end{proof}

\begin{my_lemma}
Let $\Gamma$ be an~aspherical graph and let $D$ be a~minimal non-trivial spherical diagram over the~presentation $\langle S\mid R_f \rangle$. Then every face of~$\bar D$ is simply connected and the~boundary of~every face of~$\bar D$ lifts to a~simple closed path in~the~graph $\Gamma$.
\end{my_lemma}

\begin{proof}[Proof]
Assume the~contrary. Then there are two types of~bad faces: not simply connected faces and simply connected faces with a~not simple boundary path which does not lift to a~simple closed path in~the~graph. Note that since a~not originating skeleton does not have spurs every face with a~not simple boundary encloses some subdiagram which has at least one face. Similarly every not simply connected face does.

Consider a~face $\Pi$ which is the~innermost bad face. It means that a~subdiagram $\Delta$ enclosed by this face contains only simply connected faces which boundaries lift to simple closed paths.

First examine the~case when $\Pi$ is a~simple connected face. As before let $p$ be a~subpath of~$\partial \Pi$ which lifts to a~simple closed path and let $s$ and $t$ be respectively the~start and the~end of~the~path $p$. If $s \neq t$ then acting as in~the~proof of~Lemma $4$ we obtain a~contradiction with minimality of~$D$. Thus $s=t$.

Let $q$ be such a~subpath of~$\partial \Pi$ that $pq \in \partial \Pi^{+}$. Note that either $p$ or $q$ encloses some subdiagram $\Delta'$ of~the~diagram $\Delta$. If $q$ does then we can again assume that there is a~subpath $p'$ in~$q$ with the~end points $s'$ and $t'$ which lifts to a~simple closed path. Arguing as before we obtain $s'=t'$. And now $p'$ encloses some subdiagram $\Delta'$ of~the~diagram $\Delta$. Note that $\Delta'$ as a~subdiagram of~$\Delta$ does not contains bad faces. So in~the~both cases there is a~subpath which encloses some subdiagram without bad faces and which lifts to a~simple closed path. We denote this path by $p$ and this subdiagram by $\Delta$.

Since $\Delta$ does not contain bad faces boundaries of~all faces of~$\Delta$ lift to a~simple closed paths. But the~boundary of~$\Delta$ itself, which equals to a~path $p$, lifts to a~simple closed path. Thus we can think of~$\Delta$ as a~spherical diagram over the~presentation $\langle S\mid R_s \rangle$ no edge of~which originates from the~graph. But it contradicts asphericity of~the~graph $\Gamma$.

In the~case of~a~not simple connected face we act similarly. Let $\Delta$ be a~subdiagram enclosed by $\Pi$. Due to the~choice of~$\Pi$ all faces of~$\Delta$ lift to simple closed paths. If $\partial \Delta$ lifts to a~simple closed path in~the~graph then we again contradict asphericity of~the~graph. Otherwise we take a~subpath $p$ of~$\partial \Delta$ with the~end points $s$ and $t$ which lifts to a~simple closed path in~the~graph. If $s \neq t$ then we contradict minimality of~$D$. If $s=t$ then a~subdiagram enclosed by $p$ gives a~spherical diagram over the~presentation $\langle S\mid R_s \rangle$ where no edge originates from the~graph that contradicts asphericity of~the~graph.
\end{proof}

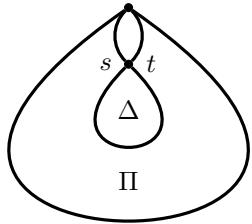
\begin{figure}[!hb]
\centering
\begin{tikzpicture}
[node_style/.style={shape=circle,draw,fill=black!100,inner sep=1pt, minimum size=3pt}]

\node (A) at (0, 0) [node_style] {};
\node (B) at (0, -0.75) [node_style] {};
\node () at (0.3, -0.75) {$t$};
\node () at (-0.3, -0.75) {$s$};
\path
    (A) edge[-, very thick, out=-45, in=45] (B)
    (A) edge[-, very thick, out=-135, in=135] (B)
	(B) edge[-, in=-135, out=-45, looseness=0.9, distance=2cm, very thick] (B)
	(A) edge[-, in=-145, out=-35, looseness=0.5, distance=6.5cm, very thick] (A);
\node () at (0, -1.35) {$\Delta$};
\node () at (0, -2.3) {$\Pi$};
\end{tikzpicture}
\caption{Illustration for Lemma 5.}
\end{figure}

\pagebreak

The following lemma can be found in~\cite{Gru15}.
\begin{my_lemma}[Gruber]
Let $\Gamma$ be a~connected $C(2)$-labelled graph. Let $R$ be the~set of~cyclic reductions of~words read a~set of~free generators of~$\pi_1(\Gamma, v)$ for some $v \in \Gamma$. Let $(D, v)$ be a~simple disk diagram over $R$ with freely trivial boundary word such that every interior edge originates from $\Gamma$. Then any sequence obtained from $(D, v)$ is a~trivial identity sequence.
\end{my_lemma}

Now we can prove the~theorem.

\begin{proof}[Proof of~the~theorem]
Assume the~contrary, that there exists non-trivial identities. Let $D$ be a~minimal non-trivial spherical diagram over the~presentation $\langle S\mid R_f \rangle$. Consider a~not originating skeleton of~$D$. If the~skeleton is empty then all edges of~the~diagram originate from the~graph. In~this case we can unglue $D$ along some edge connecting distinct vertices and we obtain a~simple disk diagram with a~trivial boundary label all edges of~which originate from the~graph $\Gamma$. Clearly, all these edges lift to the~same connected component of~the~graph. Then, by the~previous lemma, $D$ is a~diagram for a~trivial identity that contradicts its definition. Thus the~not originating skeleton is not empty. By the~all previous lemmas every face of~$\bar D$ is simply connected and lifts to a~simple closed path of~the~graph. Therefore $\bar D$ is a~spherical diagram over the~presentation $\langle S\mid R_s \rangle$ where no edge originates from the~graph that contradicts asphericity of~the~graph.
\end{proof}

\section{Corollary of~the~main result}
\label{corollary}

\subsection{Small cancellation conditions}

Recall that classical small cancellations conditions $C(q)\&T(p)$ (see, for example, \cite{LiShu77}) mean that every face in any reduced spherical diagram consists of at least $q$ arcs and every vertex has degree at least $p$ or equal to $2$. Graphical analogue of condition $C(q)$ was already introduced in the article, so it only remains to formulate an analogue of condition $T(p)$.

Let $\Gamma$ be a labelled by the set $S$ $C(2)$-graph. Let $r_1$ and $r_2$ be two elements of $R_s$. We say that $r_1$ and $r_2$ mutually originate from the graph if $r_1 = r_1'c$, $r_2 = c^{-1} r_2'$ and lifts of $c$ in the graph via $r_1$ and via $r_2$ coincide (recall that $r_1$ and $r_2$ have unique lifts since they are elements of $R_s$ and $\Gamma$ is a $C(2)$-graph). Now we can modify a classical definition from \cite{LiShu77} to obtain a definition of graphical condition $T(q)$.

\begin{my_def}[Graphical condition $T(p)$]
	Let $\Gamma$ be a labelled by the set $S$ $C(2)$-graph and let ${3 \leqslant h < p}$. Assume $r_1, \ldots, r_h$ to be elements of $R_s$ such that consecutive elements $r_i$, $r_{i+1}$ are not mutually originating from the graph. Then at least one product $r_1 r_2,\,\ldots ,\,r_{h-1}r_h,\,r_h r_1$ is reduced.
\end{my_def}

This definition preserves geometric meaning of condition $T(p)$: every inner vertex of any graphically reduced diagram has degree at least $p$ or equal to $2$.

Now we show that each graphical condition $C(6)[\&T(3)]$, $C(4)\&T(4)$ or $C(3)\&T(6)$ implies asphericity of a graph. To do this we will need a~notion of~$[p, q]$-diagrams introduced in~\cite{LiShu77}.

Let $p$ and $q$ be positive integers such that $1/p\,+\,1/q = 1/2$. Degree of~a~face is a~number of~edges in~its boundary path. A~face of~a~diagram is called interior if it has no common edges with the~boundary of~the~diagram. If $D$ is a~non-empty diagram such that each interior vertex of~$D$ has degree at least $p$ and all faces of~$D$ have degree at least $q$, then $D$ is called $[p, q]$-diagram. If $D$ is a~non-empty diagram such that each interior vertex of~$D$ has degree at least $p$ and each interior face of~$D$ have degree at least $q$, then $D$ is called $(p, q)$-diagram.   

Let us show that there exists no simple spherical $(p, q)$-diagram. For the number $c$ of faces and the number $d$ of edges in such a diagram we have an inequality $c \leqslant 2 d/q$. On the other hand, for the number $v$ of vertices and the number $d$ of edges in such a diagram we have an inequality $v \leqslant 2 d/p$. Summing these inequalities and recalling that $1/p\,+\,1/q = 1/2$, we obtain a contradiction with Euler's formula.

We also need an~operation of~“forgetting” vertices of~degree $2$. Let $D$ be a~diagram without spurs. We call an arc in the diagram $D$ full if its endpoints have degree more than $2$. Then if we replace every full arc in~$D$ by an~edge, we obtain a~diagram where no vertex has degree equal to $2$. Now we are ready to prove that each graphical condition $C(6)[\&T(3)]$, $C(4)\&T(4)$ or $C(3)\&T(6)$ implies asphericity of~the~graph.

\begin{proof}[Proof]
Assume the~contrary, that there exists a~simple spherical diagram $D$ over the~presenta\-tion $\langle S\mid R_s \rangle$ no edge of~which originates from the~graph. By Lemma $1$, the~diagram $D$ has no spurs. Note that for any two adjacent faces $\Pi_1$ and $\Pi_2$ their common boundary is a~piece because it does not originate from the~graph. Due to graphical $C(q)$-condition it means that every face consists of~at least $q$ arcs. Due to graphical condition $T(p)$ every vertex of $D$ has degree at least $p$ or equal to $2$. Therefore after “forgetting” vertices of~degree $2$ the~diagram $D$ becomes a~simple spherical $(p, q)$-diagram which does not exist. A~contradiction.
\end{proof}

\begin{my_corollary}
If $\Gamma$ satisfies any of the conditions $C(6)$, $C(4)\&T(4)$ or $C(3)\&T(6)$, then the~presentation $\langle S\mid R_f \rangle$ is aspherical.
\end{my_corollary}

For condition $C(6)$ this result was obtained for the first time by Dominic Gruber in~\cite{Gru15}.

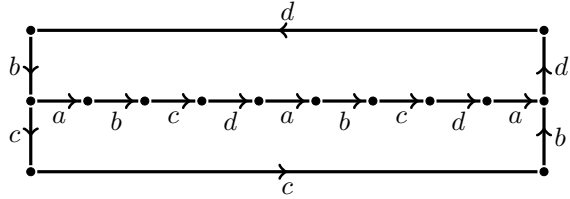
\begin{figure}[!t]
    \centering

    \begin{tikzpicture}[scale=0.75]
    \node[inner sep=0pt, minimum size=4.5pt] (A0) at (0, 0) {};
    \node[inner sep=0pt, minimum size=4.5pt] (A1) at (1, 0) {};
    \node[inner sep=0pt, minimum size=4.5pt] (A2) at (2, 0) {};
    \node[inner sep=0pt, minimum size=4.5pt] (A3) at (3, 0) {};
    \node[inner sep=0pt, minimum size=4.5pt] (A4) at (4, 0) {};
    \node[inner sep=0pt, minimum size=4.5pt] (A5) at (5, 0) {};
    \node[inner sep=0pt, minimum size=4.5pt] (A6) at (6, 0) {};
    \node[inner sep=0pt, minimum size=4.5pt] (A7) at (7, 0) {};
    \node[inner sep=0pt, minimum size=4.5pt] (A8) at (8, 0) {};
    \node[inner sep=0pt, minimum size=4.5pt] (A9) at (9, 0) {};

    \node[inner sep=0pt, minimum size=4.5pt] (B1) at (0, 1.25) {};
    \node[inner sep=0pt, minimum size=4.5pt] (B2) at (9, 1.25) {};

    \node[inner sep=0pt, minimum size=4.5pt] (C1) at (0, -1.25) {};
    \node[inner sep=0pt, minimum size=4.5pt] (C2) at (9, -1.25) {};

    \filldraw (A0) circle (2pt);
    \filldraw (A1) circle (2pt);
    \filldraw (A2) circle (2pt);
    \filldraw (A3) circle (2pt);
    \filldraw (A4) circle (2pt);
    \filldraw (A5) circle (2pt);
    \filldraw (A6) circle (2pt);
    \filldraw (A7) circle (2pt);
    \filldraw (A8) circle (2pt);
    \filldraw (A9) circle (2pt);
    \filldraw (B1) circle (2pt);
    \filldraw (B2) circle (2pt);
    \filldraw (C1) circle (2pt);
    \filldraw (C2) circle (2pt);

    \draw[->-=0.9, very thick] (A0) -- node[below] {$a$} (A1);
    \draw[->-=0.9, very thick] (A1) -- node[below] {$b$} (A2);
    \draw[->-=0.9, very thick] (A2) -- node[below] {$c$} (A3);
    \draw[->-=0.9, very thick] (A3) -- node[below] {$d$} (A4);
    \draw[->-=0.9, very thick] (A4) -- node[below] {$a$} (A5);
    \draw[->-=0.9, very thick] (A5) -- node[below] {$b$} (A6);
    \draw[->-=0.9, very thick] (A6) -- node[below] {$c$} (A7);
    \draw[->-=0.9, very thick] (A7) -- node[below] {$d$} (A8);
    \draw[->-=0.9, very thick] (A8) -- node[below] {$a$} (A9);

    \draw[-<-=0.5, very thick] (A0) -- node[left] {$b$} (B1);
    \draw[->-=0.5, very thick] (A0) -- node[left] {$c$} (C1);

    \draw[->-=0.5, very thick] (A9) -- node[right] {$d$} (B2);
    \draw[-<-=0.5, very thick] (A9) -- node[right] {$b$} (C2);

    \draw[-<-=0.5, very thick] (B1) -- node[above] {$d$} (B2);
    \draw[->-=0.5, very thick] (C1) -- node[below] {$c$} (C2);
    \end{tikzpicture}
    
    \caption{An example of $C(4)\&T(4)$-graph}
    \label{example_main:fig}
\end{figure}

Consider an example. Let $\Gamma$ be a graph as in figure \ref{example_main:fig}., which defines a group $\langle a, b, c, d \mid (abcd)^2ad^2b = 1, (abcd)^2a = c^2b \rangle$. It is easy to check that this presentation doesn't satisfy classical condition $C(6)$. Neither it satisfies classical condition $T(4)$, because we can consider relations $c^2ba^{-1}(abcd)^{-2}$, $(abcd)^2ad^2b$ and $b^{-1}a^{-1}c^2ba^{-1}(abcd)^{-1}d^{-1}c^{-1}$. The graph $\Gamma$ doesn't satisfy graphical condition $C(6)$ because, for example, a simple cycle $(abcd)^2ad^2b$ can be written as $abcda \cdot bcda \cdot d \cdot d \cdot b$.

Nevertheless, we can check that $\Gamma$ satisfies graphical condition $C(4)\&T(4)$ and therefore this presentation is aspherical due to previous corollary. It is easy to check graphical condition $C(4)$.

Let us check that $\Gamma$ satisfies graphical condition $T(4)$. We should show that for any triple of relations $r_1, r_2, r_3 \in R_s$, such that each pair of consecutive elements have a reduction, a mutually originating pair exists. To check that each pair in triple $r_1, r_2, r_3$ has a reduction it is sufficient to consider only first and last letters of each relation. And since all cyclic shifts lie in $R_s$ it is sufficient to consider only two-letter subwords of relations. 

Note that any two-letter subword of a relation from $R_s$ looks like $xy$, where either $x, y \in \{a, b, c, d\}$ or $x, y \in \{a^{-1}, b^{-1}, c^{-1}, d^{-1}\}$, only except subwords $ab^{-1}$, $ba^{-1}$, $a^{-1}c$, $c^{-1}a$. It is clear that if in a triple of subwords each pair have a reduction then at least one of them should be an exceptional subword. Considering possible cases we obtain only $4$ appropriate triples: $(ab^{-1}, bd, d^{-1}a^{-1}), (a^{-1}c$, $c^{-1}b^{-1}, ba)$ and inverse to them. Finally note that a letter $a$ originates in all this triples.

\subsection{Car-crash lemma}
Shortly consider one more method for proving diagrammatic asphericity, which is based on a topological lemma from $\cite{Kl93}$. Let $S$ be a simple spherical diagram and let there be a moving point (a car) on the boundary of some face of the diagram. We say that a car moves properly if it moves along the boundary in the positive direction continuously, perpetually, with no stops, no reverses and visiting every point of the boundary infinitely many times. 

\begin{my_lemma}[Klyachko]
Let $S$ be a simple spherical diagram. Let there be a car on the boundary of each face and let the cars move properly. Then there exists at least $2$ points of the sphere in which complete collision happens. A collision is called complete, if in a point of multiplicity $k$ collide $k$ cars simultaneously. 
\end{my_lemma} 

This lemma can be used to prove diagrammatic asphericity in the following way. Assume that a motion is defined for each relation such that for any diagram over that presentation collisions occurs only on common boundaries of reducible pairs of faces. Then there exist no reduced spherical diagrams because otherwise there exists a motion of cars on a sphere without collisions that contradicts the lemma.

Consider, for example, a presentation $\langle a, b \mid aba^{-1}b^{-1} \rangle$. Let us prove that it is aspherical. Define a motion of cars. Let each car evenly moves along its face with period $1$. If a car moves along a face with relation $aba^{-1}b^{-1}$ then first it traverses an edge with a letter $a$, then it traverses an edge with a letter $b$ and so on. If a car moves along a face with relation $bab^{-1}a^{-1}$ then first it traverses an edge with a letter $a$, then it traverses an edge with a letter $b^{-1}$ and so on.

In the moments of time from $(k, k + 1/4)$ all cars move along edges with a letter $a$. In the moments of time from $(k + 2/4, k + 3/4)$ all cars move along edges with a letter $a^{-1}$. Hence there can be no collisions in these moments of time. In the moments of time from $(k + 1/4, k + 2/4)\cup(k + 3/4, k + 1)$ all cars move along edges with letters $b^{\pm1}$. Hence there can be no complete collision in a vertex because if it occurs in the moment $t$ then in the moment $t + \varepsilon$ one part of the cars should move along edges with letters $a^{\pm1}$ and another part should move along edges with letters $b^{\pm1}$ that is impossible. In the moments of time $k + 3/8$ and $k + 7/8$ a collision on an edge with a letter $b$ can occur but in this case a collision occurs between a reducible pair of faces. Thus collisions occurs only on common boundaries of reducible pairs of faces that imply asphericity of the presentation $\langle a, b \mid aba^{-1}b^{-1} \rangle$.

This method can be transfered to the graphical case in the following way: let us allow collisions to occur on edges originating from the graph.

\begin{my_corollary}
    Let $\Gamma$ be a $C(2)$-graph and let a proper motion be given for each (up to conjugation) element of $R_s$ such that complete collisions in diagrams over $\langle S\mid R_s \rangle$ occur only on edges originating from the graph and on vertices incident to originating edges. Then the presentation $\langle S\mid R_f \rangle$ is aspherical.
\end{my_corollary}

\renewcommand{\refname}{References}

\end{document}